\documentclass[hidelinks,onefignum,onetabnum]{siamart251216}
\usepackage{amsmath,amssymb,mathtools}
\renewcommand{\boxed}[1]{{\displaystyle #1}}
\usepackage{booktabs,tabularx,array,multirow}
\usepackage{graphicx,algorithm,algpseudocode}
\usepackage{float,placeins}
\usepackage{chngcntr}
\counterwithout{equation}{section}
\graphicspath{{figures/}}

\headers{Asynchronous Concurrent Algebraic Multigrid}{C. Feng, S. Liu, L. Zhang, and Z. Zhu}
\title{A Parallel Algebraic Multigrid Preconditioned Iterative Method Based on an MPI Process Grouping Strategy for Asynchronous Concurrent Execution of the Setup and Solve Phases}
\author{Chunsheng Feng\thanks{School of Mathematics and Computational Science, Xiangtan University, Xiangtan 411105, Hunan, China; National Center for Applied Mathematics in Hunan, Hunan Shaofeng Institute for Applied Mathematics, Hunan Research Center of the Basic Discipline Fundamental Algorithmic Theory and Novel Computational Methods, Xiangtan 411105, China. \email{spring@xtu.edu.cn}.}
\and Shuai Liu\thanks{School of Mathematics and Computational Science, Xiangtan University, Xiangtan 411105, Hunan, China. \email{202631510179@smail.xtu.edu.cn}.}
\and Linjie Zhang\thanks{School of Mathematics and Computational Science, Xiangtan University, Xiangtan 411105, Hunan, China. \email{202231510136@smail.xtu.edu.cn}.}
\and Zirui Zhu\thanks{School of Mathematics and Computational Science, Xiangtan University, Xiangtan 411105, Hunan, China. Corresponding author. \email{202521511281@smail.xtu.edu.cn}.}}

\begin{document}
\raggedbottom
\maketitle
\begin{abstract}
Algebraic multigrid (AMG) is one of the most effective algorithms for solving large sparse linear systems and consists of a setup phase followed by a solve phase. As supercomputer performance improves and the number of MPI processes increases, AMG faces serious challenges in algorithmic scalability and parallel scalability. To address these issues, we propose an asynchronous concurrent algebraic multigrid (ACAMG) model based on an MPI process grouping strategy for setup and solve. The setup and solve phases execute asynchronously and concurrently on separate MPI process groups. This organization can both preserve AMG hierarchy quality through an appropriate choice of the setup-group size and achieve time overlap between the two phases through asynchronous concurrent execution. Based on this model, we design and implement the MPI ACAMG algorithm. Theoretical analysis and numerical experiments show that the algorithm exhibits better algorithmic and parallel scalability. For example, for a three-dimensional diffusion problem with spatially localized anisotropy, ACAMG not only reduces the solve-phase iteration count compared with the conventional method, but also achieves speedups of 4.74 and 19.36 at 64 and 128 processes, respectively.
\end{abstract}

\begin{keywords}
asynchronous concurrent execution, performance model, parallel computing, algebraic multigrid
\end{keywords}
\begin{MSCcodes}
65N55, 65F08, 65F10, 65Y05
\end{MSCcodes}

\section{Introduction}

Discretization of partial differential equations by finite difference, finite element, and finite volume methods usually leads to large sparse linear systems of the form
\begin{equation}
 Ax=b,\qquad A\in\mathbb{R}^{n\times n}.
 \label{eq:linear_system}
\end{equation}
Multigrid methods attenuate high-frequency error through relaxation and treat low-frequency error through coarse-grid correction, making them among the most effective methods for solving such discretized systems. They are commonly classified as geometric multigrid (GMG) or algebraic multigrid (AMG) \cite{brandt1977,ruge1987}. GMG relies on a geometric grid hierarchy, whereas AMG does not require an explicitly supplied geometric hierarchy and is therefore well suited to complex geometries and unstructured discretizations. AMG is an important research direction in high-performance scientific computing and remains an active and challenging topic.

Over time, several representative approaches to AMG coarse-space construction have been developed. Classical AMG commonly employs Ruge--St\"uben (RS) coarsening with C/F splitting based on strong connections; parallel coarsening strategies include CLJP, Falgout, PMIS, and HMIS \cite{ruge1987,cleary1998,henson2002,moxu2007,desterck2006}. Adaptive AMG uses test vectors or error components that are difficult to eliminate to identify low-energy modes to be retained in the coarse space; representative methods include Adaptive AMG and Bootstrap AMG \cite{chartier2003spectral,brezina2005,brezina2006,brandt2011bootstrap}. Reduction-based AMG (AMGr) constructs coarse-grid correction using ideal interpolation, ideal restriction, and block elimination \cite{xu2017amg,alsalti2024amgr}. Aggregation-based AMG groups neighboring degrees of freedom into aggregates and constructs the coarse space from them, giving rise to unsmoothed aggregation (UA) and smoothed aggregation (SA) methods \cite{vanek1996,blatt2013massively}.

With the development of high-performance computing architectures and parallel computing techniques, parallel algebraic multigrid methods have been widely studied. AMG-preconditioned Krylov subspace methods have become a mainstream approach for solving large sparse linear systems. Research on solution efficiency at large parallel scales mainly addresses two aspects: algorithmic scalability, namely whether iteration counts remain stable as the parallel scale increases, and parallel scalability, namely whether the AMG setup phase maintains strong scaling as computational resources increase.

To improve AMG algorithmic scalability, existing studies mainly refine parallel coarsening, interpolation, and relaxation to enhance coarse-space quality and reduce the effects of increasing parallel scale on preconditioner quality. Cleary et al. and Xu and Mo improve parallel coarsening quality by coordinating coarse-point selection near process boundaries \cite{cleary1998,xu2017interface}, while Adams et al. and D'Ambra et al. study polynomial relaxation suitable for parallel computing \cite{adams2003polynomial,dambra2025polynomial}. De Sterck et al. and Yang introduce distance-two and long-range interpolation for aggressive coarsening \cite{desterck2008distance,yang2010longrange}; D'Ambra et al. and Janna et al. enhance coarse-space approximation by improving aggregate quality and using energy-minimizing interpolation, respectively \cite{dambra2023hybrid,janna2023energy}. Nevertheless, iteration counts can still increase with parallel scale in large-scale strong-scaling experiments. For a three-temperature radiation hydrodynamics problem with a fixed size of approximately $6.30\times10^6$ degrees of freedom, Zong et al. observed that the iteration count of a BoomerAMG-preconditioned Krylov method increased from 58 to 106 as the core count increased from 128 to 2048 \cite{zong2025structmg}. This indicates that, for a fixed problem size, AMG algorithmic scalability may still deteriorate as the parallel scale increases.

To improve AMG parallel scalability, existing studies mainly reduce parallel overhead through coarse-level computational organization and communication optimization. Blatt et al. progressively aggregate coarse-grid problems onto fewer active processes to increase coarse-level computational granularity \cite{blatt2013massively}; Falgout and Schroder use non-Galerkin coarse-grid operators to reduce coarse-level operator complexity and communication overhead \cite{falgout2014nongalerkin}; and Semi-StructMG reduces hierarchy complexity through multidimensional coarsening \cite{zong2025semistructmg}. For communication optimization, additive AMG reduces communication waits by increasing concurrency across levels \cite{vassilevski2014additive}; Bienz et al. reduce inter-node messages through node-aware communication \cite{bienz2020nodeaware}; and Yuan et al. propose communication-reduced AMG (CRAMG), which reduces matrix--vector products and neighbor communication in multigrid cycles through algebraic reorganization \cite{yuan2025cramg}. Nevertheless, strong-scaling experiments still show that communication and coarse-level costs may gradually dominate total time. For a fixed-size three-dimensional Laplace problem, Bienz et al. found that, beyond approximately 512 processes, communication costs in the standard implementation offset the savings from reduced local computation, and further increases in the process count led to longer total solution times \cite{bienz2020nodeaware}. As parallel granularity decreases, AMG parallel scalability may still deteriorate, and total time may saturate or even increase.

Conventional distributed AMG--FGMRES usually completes AMG setup before starting Krylov iterations on the same MPI process group. Its total setup and solve time therefore consists of two sequential contributions:
\begin{equation}
 T_{\mathrm{total}}=T_{\mathrm{setup}}+T_{\mathrm{solve}}.
\end{equation}
As the parallel scale increases, limited setup strong scaling and deterioration of AMG preconditioner quality may increase setup time and the outer iteration count, respectively, jointly limiting end-to-end speedup.

In terms of execution organization, asynchronous iteration allows computing units to advance independently using the information already available \cite{frommer2000}. Wolfson-Pou and Chow apply this approach to multigrid solves to reduce synchronization waits \cite{wolfson2019async}. Treister and Yavneh update the hierarchy during the solve \cite{treister2011}; for sequences of linear systems, Xu et al. reduce repeated setup through adaptive preconditioning and on-demand construction \cite{xu2016adaptive,xu2020alphasetup}. These studies improve execution through changes to iteration synchronization or hierarchy updates.

In contrast to these methods, we draw on the task organization of multiple-instruction, multiple-data (MIMD) execution and propose an asynchronous concurrent algebraic multigrid (ACAMG) model based on an MPI process grouping strategy. AMG hierarchy construction and preconditioned Krylov iteration for the same linear system are assigned to independent process groups and proceed concurrently.

The main contributions are as follows:
\begin{enumerate}
  \item We propose an asynchronous concurrent algebraic multigrid (ACAMG) model for setup and solve based on an MPI process grouping strategy. We then develop the corresponding time model and analyze the speedup mechanisms and conditions for reducing total stage time.
  \item Based on the ACAMG model, we design and implement the MPI ACAMG--FGMRES algorithm. Numerical experiments verify the effectiveness of the proposed model and algorithm and demonstrate that the method can improve AMG algorithmic scalability and parallel scalability and reduce end-to-end solution time.
\end{enumerate}

The remainder of this paper is organized as follows. Section~2 presents the ACAMG model, total-time function, and speedup conditions. Based on this model, section~3 develops the MPI ACAMG algorithm using FGMRES. Section~4 reports numerical experiments, and section~5 concludes the paper.

\section{ACAMG model}
\label{sec:acamg_model}

\subsection{Model of conventional AMG-preconditioned Krylov iteration}

Consider the linear system $Ax=b$ distributed over $P$ MPI processes. A conventional AMG-preconditioned Krylov subspace method performs AMG setup followed by Krylov iteration with a fixed AMG preconditioner on the same process group:
\begin{equation}
\operatorname{setup}(P) \longrightarrow \operatorname{solve}(P).
\label{eq:v106_trad_workflow}
\end{equation}

Let $T_{\mathrm{setup}}^{\mathrm{trad}}(P)$ and $T_{\mathrm{solve}}^{\mathrm{trad}}(P)$ denote the setup and solve times, respectively. The total stage time is
\begin{equation}
T_{\mathrm{trad}}(P) = T_{\mathrm{setup}}^{\mathrm{trad}}(P) + T_{\mathrm{solve}}^{\mathrm{trad}}(P).
\label{eq:v106_trad_stages}
\end{equation}
Writing the solve time in terms of the outer iteration count and mean iteration time gives
\begin{equation}
T_{\mathrm{solve}}^{\mathrm{trad}}(P) = I_{\mathrm{trad}}(P) T_{\mathrm{iter}}^{\mathrm{trad}}(P).
\label{eq:v105_author_5}
\end{equation}
Here $I_{\mathrm{trad}}(P)$ is the number of outer iterations required by the conventional AMG-preconditioned Krylov method to satisfy the prescribed stopping criterion. The mean time per outer iteration, $T_{\mathrm{iter}}^{\mathrm{trad}}(P)$, with a fixed AMG preconditioner includes preconditioner application, matrix--vector multiplication, orthogonalization, and the associated communication costs.

Substituting into \eqref{eq:v106_trad_stages} yields
\begin{equation}
T_{\mathrm{trad}}(P) = T_{\mathrm{setup}}^{\mathrm{trad}}(P) + I_{\mathrm{trad}}(P) T_{\mathrm{iter}}^{\mathrm{trad}}(P).
\label{eq:traditional_endpoint}
\end{equation}
Changing the process count can affect setup computation and communication costs, as well as AMG preconditioning effectiveness, thereby changing the outer iteration count.

\subsection{ACAMG model and its two-stage time model}

ACAMG uses an MPI process grouping strategy to assign AMG construction and Krylov iteration to separate process groups. Let $P$ be the total number of MPI processes, and let the setup-group size $M$ be the model variable, with
\begin{equation}
M\in\{0,1,\ldots,P-1\}.
\label{eq:unified_m_domain}
\end{equation}

For $M=0$, all $P$ processes perform AMG construction followed by Krylov iteration, corresponding to the conventional model in section~2.1. For $1\le M\le P-1$, the processes are partitioned into a setup group and a solve group containing $M$ and $P_s=P-M$ processes, respectively, so that
\begin{equation}
P=M+P_s.
\label{eq:v106_resource_sum}
\end{equation}
The setup group constructs the AMG hierarchy, while the solve group performs preconditioned Krylov iterations; the two groups execute concurrently, that is,
\begin{equation}
\operatorname{setup}(M)\parallel\operatorname{solve}(P_s).
\label{eq:setup_solve_concurrency}
\end{equation}
The setup group transfers the progressively constructed hierarchy to the solve group, allowing the solve group to update its preconditioner using the currently available hierarchy. The algorithm can terminate early if the solve group satisfies the convergence criterion. If the complete hierarchy becomes available first and the solve group has not converged, the complete hierarchy and the current iterate are redistributed to all $P$ processes for continued preconditioned Krylov iteration. Accordingly, ACAMG consists of an asynchronous concurrent stage and a post-asynchronous all-process solve stage.

Take the start of concurrent execution by the two groups as the time origin. The AMG levels are numbered $0,1,\ldots,L$ from finest to coarsest; level~0 is the initial fine level and level~$L$ is the coarsest level. Let $\Delta T_\ell^{\mathrm{build}}(M)$ denote the elapsed time between completion of levels $\ell-1$ and $\ell$. The construction completion time for level~$\ell$ is then
\begin{equation}
\boxed{
t_\ell^{\mathrm{build}}(M)
=
\sum_{k=1}^{\ell}
\Delta T_k^{\mathrm{build}}(M),
\qquad 1\leq\ell\leq L.
}
\label{eq:level_build_time}
\end{equation}
The initial fine-grid operator is given before the asynchronous stage begins, so $t_0^{\mathrm{build}}(M)=0$.

After constructing each level, the setup group transfers the corresponding data to the solve group through a data-flow mechanism. For level~$\ell$, let $T_\ell^{\mathrm{trans}}(M,P_s)$ denote the time from completion of construction to completion of transfer, and let $T_\ell^{\mathrm{adopt}}(M,P_s)$ denote the additional time required for assembly, communication preparation, and group-wide confirmation before the level can safely be used in a preconditioner application. If some assembly overlaps with transfer, adoption time includes only the work remaining after transfer is complete. The data-flow time for level~$\ell$ is defined by
\begin{equation}
\boxed{
T_\ell^{\mathrm{flow}}(M,P_s)
=
T_\ell^{\mathrm{trans}}(M,P_s)
+
T_\ell^{\mathrm{adopt}}(M,P_s).
}
\label{eq:level_flow_components}
\end{equation}
Thus, the time at which level~$\ell$ becomes available to the solve group is
\begin{equation}
\boxed{
\begin{aligned}
t_\ell^{\mathrm{ready}}(M)
&=
t_\ell^{\mathrm{build}}(M)
+
T_\ell^{\mathrm{flow}}(M,P_s)\\
&=
\sum_{k=1}^{\ell}
\Delta T_k^{\mathrm{build}}(M)
+
T_\ell^{\mathrm{flow}}(M,P_s),
\qquad 1\leq\ell\leq L.
\end{aligned}
}
\label{eq:level_flow_ready}
\end{equation}
The solve group collectively adopts a new level before a preconditioner application. If data arrive during the current iteration, adoption waits until that iteration has finished. This waiting time is included in $T_\ell^{\mathrm{adopt}}$, so $t_\ell^{\mathrm{ready}}$ is the time at which the level can actually be used collectively in a preconditioner application.

After constructing the complete AMG hierarchy, the setup group transfers the final level through the data-flow mechanism and sends a construction-completion notification. The solve group obtains the complete hierarchy after adopting level~$L$ and confirming this notification. Let $T_{\mathrm{notify}}(M,P_s)$ denote the additional time, after level~$L$ becomes available, required to confirm completion and coordinate a safe stage transition. The time at which the complete hierarchy is ready in the solve group is
\begin{equation}
\boxed{
\begin{aligned}
t^{\mathrm{full}}(M)
&=
t_L^{\mathrm{ready}}(M)
+
T_{\mathrm{notify}}(M,P_s)\\
&=
\sum_{k=1}^{L}
\Delta T_k^{\mathrm{build}}(M)
+
T_L^{\mathrm{flow}}(M,P_s)
+
T_{\mathrm{notify}}(M,P_s).
\end{aligned}
}
\label{eq:v106_full_ready}
\end{equation}
Here $T_{\mathrm{notify}}(M,P_s)$ includes only notification and coordination work remaining after level~$L$ becomes available; it excludes the transfer and adoption work already included in $T_L^{\mathrm{flow}}(M,P_s)$. 

For convergence during the asynchronous concurrent stage, suppose that the solve group converges while level~$\ell$ is the deepest available level, where $0\le\ell\le L$. If level~$L$ has been adopted but the construction-completion notification has not yet been confirmed, the solve group can still iterate using the current hierarchy. Let $I_\ell^{\mathrm{conv}}(M,P_s)$ be the number of Krylov iterations performed between adoption of this level and convergence, and let $\overline T_\ell^{\mathrm{iter}}(M,P_s)$ be their average time. The convergence time is
\begin{equation}
\boxed{
\begin{aligned}
t^{\mathrm{conv}}(M)
={}&
t_\ell^{\mathrm{ready}}(M)\\
&+
I_\ell^{\mathrm{conv}}(M,P_s)
\overline T_\ell^{\mathrm{iter}}(M,P_s).
\end{aligned}
}
\label{eq:async_first_convergence}
\end{equation}
We take $t_0^{\mathrm{ready}}(M)=0$. The average iteration time accounts for the actual iteration cost over the interval from adoption of the current level to convergence, including computation, communication, and synchronization.

If the solve group converges first, it must notify the setup group to stop further construction, complete outstanding transfers, and coordinate termination between the groups. Let $T_{\mathrm{finish}}(M)$ denote the time required for this work. The actual duration of the asynchronous concurrent stage is therefore
\begin{equation}
\boxed{
T_{\mathrm{async}}(M)
=
\begin{cases}
t^{\mathrm{conv}}(M)+T_{\mathrm{finish}}(M),
&
t^{\mathrm{conv}}(M)\leq t^{\mathrm{full}}(M),\\[2mm]
t^{\mathrm{full}}(M),
&
t^{\mathrm{full}}(M)<t^{\mathrm{conv}}(M).
\end{cases}
}
\label{eq:async_stage_time}
\end{equation}
If the two events occur simultaneously, the convergence branch takes precedence. If the solve group does not converge before the complete hierarchy is ready, we take $t^{\mathrm{conv}}(M)=+\infty$ in this event comparison. If early convergence stops further AMG construction, we take $t^{\mathrm{full}}(M)=+\infty$. The term $T_{\mathrm{finish}}(M)$ distinguishes communication completion from the instant at which the convergence criterion is first satisfied, so that the model accounts for the actual time needed to end the asynchronous stage.

If the complete AMG hierarchy is collectively available in the solve group and convergence has not occurred, the complete hierarchy and the iterate improved during the asynchronous stage are redistributed to all $P$ processes. FGMRES is reinitialized and continues with the complete hierarchy as its preconditioner.

Let $T_{\mathrm{flow}}^{\mathrm{switch}}(M,P_s,P)$ denote the time required for this data-flow operation, including target-distribution mapping, interprocess transfer, receive-side assembly, and the communication preparation and confirmation needed for a usable distribution. Let $T_{\mathrm{reinit}}(M)$ denote the time required to reinitialize the FGMRES state after redistribution. The stage-transition time is
\begin{equation}
\boxed{
T_{\mathrm{switch}}(M)
=
T_{\mathrm{flow}}^{\mathrm{switch}}(M,P_s,P)
+
T_{\mathrm{reinit}}(M).
}
\label{eq:v106_switch}
\end{equation}

Let $I_{\mathrm{post}}(M)$ be the remaining number of Krylov iterations in the all-process solve stage, and let $T_{\mathrm{iter}}^{\mathrm{post}}(M,P)$ be their average time. The duration of this stage is
\begin{equation}
\boxed{
T_{\mathrm{post}}(M;P)
=
\begin{cases}
0,
&
t^{\mathrm{conv}}(M)\leq t^{\mathrm{full}}(M),\\[2mm]
T_{\mathrm{switch}}(M)
+
I_{\mathrm{post}}(M)
T_{\mathrm{iter}}^{\mathrm{post}}(M,P),
&
t^{\mathrm{full}}(M)<t^{\mathrm{conv}}(M).
\end{cases}
}
\label{eq:post_continuation_branch}
\end{equation}
Here $I_{\mathrm{post}}(M)$ counts the Krylov iterations still required after entering the all-process solve stage. Since the hierarchy is constructed on $M$ processes, its structure and preconditioning effectiveness may differ from those of the conventional hierarchy constructed on all $P$ processes. Thus, $T_{\mathrm{iter}}^{\mathrm{post}}(M,P)$ need not equal the conventional average iteration time.

Combining the two stages gives the total ACAMG stage time for $1\le M\le P-1$:
\begin{equation}
\boxed{
T_{\mathrm{ACAMG}}(M;P)
=
T_{\mathrm{async}}(M)
+
T_{\mathrm{post}}(M;P).
}
\label{eq:two_stage_sum}
\end{equation}
Expanding this expression gives
\begin{equation}
\boxed{
T_{\mathrm{ACAMG}}(M;P)
=
\begin{cases}
t^{\mathrm{conv}}(M)+T_{\mathrm{finish}}(M),
&
t^{\mathrm{conv}}(M)\leq t^{\mathrm{full}}(M),\\[2mm]
\begin{aligned}
t^{\mathrm{full}}(M)
&+T_{\mathrm{flow}}^{\mathrm{switch}}(M,P_s,P)\\
&+T_{\mathrm{reinit}}(M)\\
&+I_{\mathrm{post}}(M)
T_{\mathrm{iter}}^{\mathrm{post}}(M,P),
\end{aligned}
&
t^{\mathrm{full}}(M)<t^{\mathrm{conv}}(M).
\end{cases}
}
\label{eq:two_stage_total}
\end{equation}
Together with the conventional model in section~2.1, this yields the unified time model
\begin{equation}
\boxed{
T_{\mathrm{total}}(M;P)
=
\begin{cases}
T_{\mathrm{trad}}(P),
&
M=0,\\[2mm]
T_{\mathrm{ACAMG}}(M;P),
&
1\leq M\leq P-1.
\end{cases}
}
\label{eq:unified_time_model}
\end{equation}
For a fixed total process count $P$, the setup-group size $M$ affects hierarchy construction time, intergroup data-flow cost, asynchronous-stage convergence progress, and the subsequent all-process solve cost. Section~2.3 analyzes how these factors affect time savings, speedup conditions, and the choice of setup-group size.

\subsection{Gain over the conventional method, speedup conditions, and optimal setup-group size}

The time saved by ACAMG relative to conventional AMG-preconditioned Krylov iteration is
\begin{equation}
T_{\rm gain}(M;P) = T_{\rm trad}(P) - T_{\rm ACAMG}(M;P), \qquad M>0.
\label{eq:acamg_gain_traditional}
\end{equation}
The solve group has already made progress during the asynchronous concurrent stage before the complete hierarchy is collectively available. The number of remaining iterations in the post-asynchronous stage therefore usually differs from that of the conventional method. Define
\begin{equation}
I_{\rm save}(M;P) = I_{\rm trad}(P) - I_{\rm post}(M),
\label{eq:equivalent_iteration_contribution}
\end{equation}
which is the reduction in subsequent all-process iterations relative to the conventional method.

If convergence occurs within the asynchronous concurrent stage,
\[
I_{\rm post}=0,
\]
and hence $I_{\rm save}=I_{\rm trad}$. The equivalent contribution $I_{\mathrm{save}}$ includes both progress made during the asynchronous stage and the effect of the setup-group size $M$ on subsequent preconditioning. It is not the number of iterations performed during the asynchronous stage.

In the continuation branch, both methods use their respective complete AMG hierarchies on all $P$ processes. When the hierarchy complexities and mean per-iteration computation and communication costs are similar, we use the approximation
\begin{equation}
T_{\rm iter}^{\rm post}(M,P) \approx T_{\rm iter}^{\rm trad}(P) \equiv T_{\rm iter}(P).
\label{eq:equal_step_assumption}
\end{equation}
We use \eqref{eq:equal_step_assumption} below to analyze trends in the time gain; the total stage time and its exact gain relation retain the mean iteration time of each method.

When the complete hierarchy becomes available before convergence, $T_{\rm async}=t^{\rm full}$. Subtracting \eqref{eq:two_stage_total} from \eqref{eq:traditional_endpoint} gives
\begin{equation}
T_{\rm gain}=T_{\rm setup}^{\rm trad}+I_{\rm trad}T_{\rm iter}^{\rm trad}-t^{\rm full}-T_{\rm switch}-I_{\rm post}T_{\rm iter}^{\rm post}.
\label{eq:post_gain_exact}
\end{equation}
Using \eqref{eq:equal_step_assumption} then gives
\begin{equation}
\begin{aligned}
 T_{\rm gain}(M;P)\approx{}&I_{\rm save}(M;P)T_{\rm iter}(P)
 -\left[t^{\rm full}(M)-T_{\rm setup}^{\rm trad}(P)\right]\\
 &-T_{\rm switch}(M).
\end{aligned}
\label{eq:equivalent_contribution_equal_step}
\end{equation}
Equation~\eqref{eq:equivalent_contribution_equal_step} separates the time gain into three terms: $I_{\rm save}T_{\rm iter}$ is the equivalent time saved by reducing subsequent complete-AMG iterations; $t^{\rm full}-T_{\rm setup}^{\rm trad}$ is the time difference between collective availability of the complete hierarchy and completion of conventional setup; and $T_{\rm switch}$ is the stage-transition time required for redistribution to all processes and solver-state reinitialization.

Speedup therefore does not require $t^{\rm full}<T_{\rm setup}^{\rm trad}$. A smaller setup group may deliver the complete hierarchy later, yet still reduce total time if asynchronous progress saves enough subsequent iterations to cover the delay and transition cost. Conversely, earlier availability alone does not guarantee speedup: limited progress on the solve group and the cost of transitioning to all processes can offset that benefit.

Expressing the two costs in iteration-equivalent units gives
\begin{equation}
I_{\rm cost}(M;P) = \frac{ t^{\rm full}(M) - T_{\rm setup}^{\rm trad}(P) + T_{\rm switch}(M) }{ T_{\rm iter}(P) }.
\label{eq:equivalent_iteration_cost}
\end{equation}
Since $t^{\rm full}-T_{\rm setup}^{\rm trad}$ may be negative, $I_{\rm cost}$ need not be nonnegative or an integer. Under the approximation in \eqref{eq:equal_step_assumption}, ACAMG is faster than the conventional method in the branch where the complete hierarchy becomes collectively available first if
\begin{equation}
I_{\rm save}(M;P) > I_{\rm cost}(M;P).
\label{eq:acamg_benefit_condition}
\end{equation}
Thus the reduction in subsequent iterations achieved during the asynchronous stage must exceed the equivalent iteration cost resulting from the change in the time at which the complete hierarchy becomes collectively available and the stage transition. Equality, $I_{\rm save}=I_{\rm cost}$, gives equal total stage times under this approximation.

If convergence occurs during the asynchronous stage, $T_{\rm post}=0$ and $T_{\rm ACAMG}=t^{\rm conv}+T_{\mathrm{finish}}$, giving
\begin{equation}
T_{\rm gain}(M;P) = T_{\rm trad}(P) - t^{\rm conv}(M)-T_{\mathrm{finish}}(M),
\label{eq:early_gain}
\end{equation}
with speedup whenever $t^{\rm conv}(M)+T_{\mathrm{finish}}(M)<T_{\rm trad}(P)$. If $t^{\rm conv}(M)+T_{\mathrm{finish}}(M)<T_{\rm setup}^{\rm trad}(P)$, ACAMG has already converged when the conventional method is still constructing its hierarchy. The time saving then includes the entire conventional solve and the remaining conventional setup time.

For fixed $P$, the choice of $M$ introduces a tradeoff. A smaller setup group leaves more processes for the asynchronous solve but may construct the hierarchy more slowly. A larger group can shorten construction time within its effective strong-scaling range, but reduces solve-group resources and may trigger an earlier transition to all processes. Changing $M$ can also change the C/F splitting, interpolation, hierarchy complexity, and transfer and adoption costs. These effects jointly determine $t^{\rm full}(M)$, $t^{\rm conv}(M)$, and $I_{\rm post}(M)$, so the total stage time may attain an interior minimum over the discrete feasible set.

The optimal setup-group size is defined by
\begin{equation}
M^\star \in \arg\min_{ M\in\{0,1,\ldots ,P-1\} } T_{\rm total}(M;P).
\label{eq:m_optimization}
\end{equation}
Since the unified model reduces to the conventional AMG-preconditioned Krylov method at $M=0$, \eqref{eq:m_optimization} includes conventional execution and ACAMG configurations at all nonzero setup-group sizes. If $M^\star=0$, conventional execution attains the lowest total stage time; if $M^\star>0$, asynchronous concurrent execution at that setup-group size is no slower than the conventional method. The setup-group-size scan below examines this tradeoff and the corresponding interior minimum.

\section{MPI ACAMG algorithm design}
\label{sec:acamg_algorithm}

Based on the two-stage time model in section~2, we design the MPI ACAMG algorithm using a unified data-flow algorithm, level-by-level AMG hierarchy construction, FGMRES with variable preconditioning, and continuation on all processes.

\subsection{Unified data-flow algorithm}
\label{sec:dataflow}

During asynchronous concurrent execution, ACAMG transfers the AMG data as they are constructed level by level from the setup group to the solve group. On entering the all-process solve stage, it redistributes the complete hierarchy and current iterate from the solve group to all processes. Although the data objects and process distributions differ, both operations require data partitioning, transfer, assembly, and confirmation. A unified data-flow algorithm therefore supports both operations.

Let $S$ and $D$ denote the source and target process groups, respectively, and let $X$ denote the data to be transferred, which may be a distributed sparse matrix, a vector, or a collection of data objects. The algorithm determines interprocess transfers from the source and target distributions.

Suppose source process $p$ owns the global row interval $[a_p,b_p)$ and target process $r$ requires $[c_r,d_r)$. The rows to be transferred between them are determined by the intersection $I_{pr}$:
\begin{equation}
\boxed{
I_{pr}=[a_p,b_p)\cap[c_r,d_r).
}
\label{eq:dataflow_row_intersection}
\end{equation}
If $I_{pr}\neq\varnothing$, process $p$ sends the corresponding data to process $r$; otherwise, no matrix or vector entries need to be sent. For a sparse matrix, the transmitted data comprise the row structure, column indices, and numerical values of the corresponding rows. Vector entries are partitioned by global index. Since matrices and vectors may have different row distributions, their transfer relations are determined separately. Expected contributions are determined from these distributions and the protocol's contribution counts; an empty intersection is represented only by control metadata when the protocol requires it.

Once these relations are determined, the source process extracts and packs the data and initiates nonblocking MPI sends. The receiver polls for incoming messages, while the sender retains buffers until the corresponding requests are confirmed complete. The target process collects available data from the source processes and assembles them by global index. For a sparse matrix, the local matrix and communication structures must also be built for the target distribution. A target process may receive data from several sources, so its local data are ready only after all expected data have been received and assembled.

We denote this unified operation by $\operatorname{DataFlow}(X,S,D)$. Algorithm~\ref{alg:dataflow} gives the basic procedure.

\begin{algorithm}[!htbp]
\caption{Unified data-flow algorithm (DataFlow)}
\label{alg:dataflow}
\footnotesize
\begin{algorithmic}[1]
\Require Data objects $X$, source group $S$, target group $D$, target row distribution
\Ensure Target data $X^{(D)}$ and availability status
\Statex \textbf{The source and target groups execute the following procedures concurrently}
\Procedure{DataFlow-Send}{$X,S,D$}
  \For{each $Q\in X$}
    \Statex \hspace{\algorithmicindent}\textbf{Each source process  $p\in S$ executes in parallel}
    \For{each $r\in D$}
      \State \textbf{Determine rows to transfer: }$I_{pr}\gets[a_p,b_p)\cap[c_r,d_r)$
      \If{$I_{pr}\neq\varnothing$}
        \State \textbf{Extract and pack data: }$B_{pr}\gets\operatorname{Pack}(Q_p,I_{pr})$
        \State \textbf{Initiate a nonblocking send: }$h_{pr}\gets\operatorname{MPI\_Isend}(B_{pr},r)$
      \Else
        \State Record the empty contribution using control metadata if required by the protocol
      \EndIf
    \EndFor
  \EndFor
  \State \textbf{Maintain send buffers: }retain outstanding requests and buffers; release buffers for completed requests
\EndProcedure
\Procedure{DataFlow-Adopt}{$X,S,D$}
  \For{each $Q\in X$}
    \Statex \hspace{\algorithmicindent}\textbf{Each target process  $r\in D$ executes in parallel}
    \State \textbf{Collect available data: }$\{B_{pr}\}\gets\operatorname{ReceiveAvailable}(Q,r)$
    \If{all expected contributions and required control metadata have arrived}
      \State \textbf{Assemble under the target distribution: }$Q_r^{(D)}\gets\operatorname{Assemble}(\{B_{pr}\})$
      \State \textbf{Build communication structures: }$\operatorname{BuildComm}(Q_r^{(D)})$
      \State \textbf{Confirm local availability: }$s_r(Q)\gets\operatorname{Ready}(Q_r^{(D)})$
    \Else
      \State \textbf{Mark as not yet ready: }$s_r(Q)\gets0$
    \EndIf
    \State \textbf{Confirm group-wide availability: }$s(Q)\gets\min_{r\in D}s_r(Q)$
  \EndFor
  \State \textbf{Return target data and status: }\Return $X^{(D)},\{s(Q)\}$
\EndProcedure
\end{algorithmic}
\end{algorithm}
\FloatBarrier

During the asynchronous concurrent stage, the setup group invokes this algorithm level by level to transfer the constructed level data to the solve group, with cost $T_\ell^{\mathrm{flow}}(M,P_s)$ in section~2.2. At the stage transition, the solve group invokes the algorithm again to redistribute the complete hierarchy and current iterate to all $P$ processes, with cost $T_{\mathrm{flow}}^{\mathrm{switch}}(M,P_s,P)$. The data-reception procedures and conditions for data availability are described in section~3.2.

\subsection{AMG construction and two-stage FGMRES solution}

\subsubsection{Level-by-level AMG hierarchy construction}
\label{sec:hierarchy-stream}

AMG constructs a hierarchy from fine to coarse by coarsening, interpolation, and coarse-grid operator construction. Let $A_0=A$ be the initial fine-grid operator, and let $P_\ell$ and $R_\ell$ be the interpolation and restriction operators at level~$\ell$. The next level is constructed using the Galerkin method:
\begin{equation}
\boxed{
A_{\ell+1}=R_\ell A_\ell P_\ell.
}
\label{eq:galerkin_product}
\end{equation}
For the symmetric problems considered here, $R_\ell=P_\ell^T$. Repeated construction yields an AMG hierarchy with $L+1$ levels.

To represent a hierarchy under different process distributions, define
\begin{equation}
\mathcal H_\ell^{(G)}
=
\left\{
A_0^{(G)},\ldots,A_\ell^{(G)};
P_0^{(G)},\ldots,P_{\ell-1}^{(G)}
\right\}.
\label{eq:amg_hierarchy_group_notation}
\end{equation}
Here $G$ is the number of processes storing the hierarchy, and $\mathcal H_\ell^{(G)}$ is the hierarchy through level~$\ell$. For $\ell=L$, $\mathcal H_L^{(G)}$ is the complete hierarchy; restriction is given by $R_i^{(G)}=(P_i^{(G)})^T$.

After constructing level~$\ell$, the setup group obtains the corresponding interpolation and coarse-grid operators. The new data are
\begin{equation}
\boxed{
Y_\ell^{(M)}
=
\left(
P_{\ell-1}^{(M)},A_\ell^{(M)}
\right),
\qquad 1\leq\ell\leq L.
}
\label{eq:hierarchy_payload}
\end{equation}
The data-flow algorithm in section~3.1 is then invoked:
\begin{equation}
\boxed{
\operatorname{DataFlow}
\left(
Y_\ell^{(M)},\mathcal P_c,\mathcal P_s
\right).
}
\label{eq:dataflow_level_call}
\end{equation}
Here $\mathcal P_c$ and $\mathcal P_s$ denote the setup and solve groups. This operation sends the new data to the corresponding solve processes according to their target distribution.

After initiating nonblocking transfer of the current level, the setup group can continue construction of the next level without waiting for assembly in the solve group. The solve group receives and assembles the new data, adding them to the hierarchy once the operators and communication structures are confirmed to be usable.

When AMG hierarchy construction is complete, the setup group sends the final level and a construction-completion notification. After adopting the final level and confirming completion, the solve group has the complete hierarchy $\mathcal H_L^{(P_s)}$. This corresponds to $t^{\mathrm{full}}(M)$ in section~2.2.

\subsubsection{FGMRES with variable AMG preconditioning and continuation on all processes}
\label{sec:fgmres-dynamic}
\label{sec:allrank-continuation}

During the asynchronous concurrent stage, the solve group uses the continually updated AMG hierarchy as its preconditioner. We therefore use FGMRES, which allows variable preconditioning \cite{saad1993,simoncini2003}.

Initially, no coarse level is available, and fine-grid smoothing provides the preconditioner application. As new levels are adopted, a V-cycle is performed using the available hierarchy. The equation on the temporary coarsest level of an incomplete hierarchy is solved approximately. The implementation supports a bounded number of GMRES iterations for the temporary coarsest-level equation; this is the configuration used in section~4.

Solve processes may receive and assemble data at different rates. A common available level must therefore be determined to ensure consistent preconditioner applications. Let $\ell_r$ be the deepest level whose data have been assembled and whose communication structures are ready on process~$r$. The common level is
\begin{equation}
\boxed{
\ell_{\mathrm{use}}
=
\min_{r\in\mathcal P_s}\ell_r.
}
\label{eq:solve_group_safe_depth}
\end{equation}
This value is obtained by an \texttt{MPI\_Allreduce} within the solve group. Before each preconditioner application, the solve group checks for and adopts new levels, then selects the preconditioner according to the deepest common level. The preconditioner for iteration~$j$ is
\begin{equation}
\boxed{
B_j=
\mathcal B\left(\mathcal H_{\ell_j}^{(P_s)}\right).
}
\label{eq:variable_amg_pc}
\end{equation}
Here $\ell_j$ is the deepest common level at the start of the $j$th preconditioner application. Once selected, the preconditioner remains unchanged during that application; subsequently adopted levels affect only later applications.

FGMRES stores the preconditioned vectors from each iteration, allowing the preconditioner to vary between steps. For the orthonormal Arnoldi vector $v_j$,
\begin{equation}
\boxed{
z_j=B_jv_j,\qquad
w_j=A_0^{(P_s)}z_j.
}
\label{eq:flex_step}
\end{equation}
Define
\begin{equation}
Z_m=[z_1,z_2,\ldots,z_m],
\qquad
V_{m+1}=[v_1,v_2,\ldots,v_{m+1}].
\label{eq:fgmres_bases}
\end{equation}
After $m$ Arnoldi iterations,
\begin{equation}
\boxed{
A_0^{(P_s)}Z_m
=
V_{m+1}\overline H_m,
}
\label{eq:flex_arnoldi}
\end{equation}
where $\overline H_m\in\mathbb R^{(m+1)\times m}$ is the upper Hessenberg matrix produced by Arnoldi orthogonalization. The approximate solution is
\begin{equation}
\boxed{
x_m=x_0+Z_my_m.
}
\label{eq:flex_update}
\end{equation}
Here $y_m=\arg\min_y\lVert\beta_0e_1-\overline H_my\rVert_2$, $\beta_0=\lVert b-A_0^{(P_s)}x_0\rVert_2$, and $e_1$ is the first coordinate vector \cite{saad1993,simoncini2003}.

Because FGMRES updates the iterate using the vectors $z_j$ actually obtained at each step, its Arnoldi relation is retained even when different iterations use different AMG hierarchies. 

After finishing the current iteration, the solve group checks convergence, collective availability of the final level, and confirmation of the construction-completion notification. If convergence has been reached, the converged solution takes precedence; otherwise, continuation on all processes begins once both readiness conditions hold.

If continuation on all processes is required, the data-flow algorithm in section~3.1 redistributes the complete hierarchy and current iterate to all $P$ processes:
\begin{equation}
\boxed{
\operatorname{DataFlow}
\left(
\left\{
\mathcal H_L^{(P_s)},x_{\mathrm{async}}^{(P_s)}
\right\},
\mathcal P_s,\mathcal P
\right).
}
\label{eq:allrank_redistribution}
\end{equation}
Here $\mathcal P$ is the set of all MPI processes. After redistribution, the resulting objects are
\begin{equation}
\boxed{
\mathcal H_L^{(P)},
\qquad
x_{\mathrm{async}}^{(P)}.
}
\label{eq:allrank_objects}
\end{equation}
This invocation uses the same mapping principles as level-by-level transfer. Its objects are now the complete hierarchy and current iterate, and continuation requires all target data to be ready. Its cost is $T_{\mathrm{flow}}^{\mathrm{switch}}(M,P_s,P)$ in section~2.2.

FGMRES is reinitialized on all $P$ processes with $x_{\mathrm{async}}^{(P)}$ as the new initial guess. The initial residual is recomputed under the new data distribution:
\begin{equation}
\boxed{
r_0^{(P)}
=
b^{(P)}
-
A_0^{(P)}x_{\mathrm{async}}^{(P)}.
}
\label{eq:allrank_restart}
\end{equation}
Initialize FGMRES from the recomputed residual and, when further iteration is required, take
\begin{equation}
v_1^{(P)}
=
\frac{r_0^{(P)}}{\|r_0^{(P)}\|_2}
\label{eq:allrank_start_vector}
\end{equation}
as the initial vector for a new Arnoldi iteration, thereby rebuilding the Krylov subspace.

FGMRES then continues on all $P$ processes with the complete AMG hierarchy $\mathcal H_L^{(P)}$ as its preconditioner until the prescribed convergence criterion is satisfied.

\FloatBarrier
\subsection{Overall algorithm}

Algorithm~\ref{alg:acamg} summarizes the construction and solution procedures in section~3.2 and invokes DataFlow-Send and DataFlow-Adopt from algorithm~\ref{alg:dataflow} for intergroup data transfer and adoption. Its two termination branches, convergence during the asynchronous stage and continuation on all processes, correspond to the two-stage time model in section~2.

\begingroup
\linespread{0.94}\selectfont
\begin{algorithm}[H]
\caption{Asynchronous concurrent algebraic multigrid (MPI ACAMG)}
\label{alg:acamg}
\footnotesize
\begin{algorithmic}[1]
\Require $(A_0,b,x_0)$, $\mathcal P$, $M\in\{0,\ldots,P-1\}$, and $\varepsilon$
\Ensure $x$ with $\eta(x)\le\varepsilon$
\If{$M=0$}
  \State \textbf{Conventional baseline:} construct the complete hierarchy on all $P$ processes
  \State \textbf{Conventional solve:} $x\gets\operatorname{FGMRES}(A_0,b,x_0;\mathcal B(\mathcal H_{L}^{(P)}),P)$
  \State \Return $x$
\Else
  \State \textbf{Partition the processes:} $\mathcal P_{\mathrm c}\mathbin{\dot\cup}\mathcal P_{\mathrm s}\gets\mathcal P$, $(|\mathcal P_{\mathrm c}|,|\mathcal P_{\mathrm s}|)=(M,P-M)$
  \State \textbf{Prepare group data and fine-grid hierarchy:} prepare $(A_0,b)$ under the two group distributions and initialize $\mathcal H_0^{(P_s)}$
  \State \textbf{Advance setup concurrently:} construct AMG level by level on $\mathcal P_{\mathrm c}$
  \Statex \hspace{\algorithmicindent}For each new $Y_\ell^{(M)}$, invoke \Call{DataFlow-Send}{$Y_\ell^{(M)},\mathcal P_{\mathrm c},\mathcal P_{\mathrm s}$}; send a notification when construction is complete
  \State \textbf{Initialize solve group:} $x\gets x_0$, $\eta\gets\lVert b-A_0x\rVert_2/\lVert b\rVert_2$; initialize FGMRES state
  \Statex \hspace{\algorithmicindent}\textbf{The FGMRES loop below runs on $\mathcal P_{\mathrm s}$ concurrently with setup}
  \While{$\eta>\varepsilon$}
    \State \textbf{Adopt candidate levels:} invoke \Call{DataFlow-Adopt}{$\{Y_\ell^{(M)}\},\mathcal P_{\mathrm c},\mathcal P_{\mathrm s}$}
    \State \textbf{Determine the commonly available level:} $\ell_j\gets\ell_{\mathrm{use}}$
    \State \textbf{Confirm construction completion:} check the notification that has arrived
    \If{$\ell_j=L$ and the construction-completion notification has been confirmed}
      \State \textbf{Retain the current iterate and continue on all processes:} \textbf{break}
    \EndIf
    \State \textbf{Apply preconditioner and store direction:} $z_j\gets\mathcal B(\mathcal H_{\ell_j}^{(P_s)})v_j$, $Z_j\gets[Z_{j-1},z_j]$
    \State \textbf{Extend flexible Arnoldi basis:} $w_j\gets A_0z_j$, $(V_{j+1},\bar H_j)\gets\operatorname{Arnoldi}(V_j,w_j)$
    \State \textbf{Update iterate:} $y_j\gets\arg\min_y\lVert\beta_0e_1-\bar H_jy\rVert_2$, $x\gets x_0+Z_jy_j$
    \State \textbf{Compute true relative residual:} $r_j\gets b-A_0x$, $\eta\gets\lVert r_j\rVert_2/\lVert b\rVert_2$
    \State Update FGMRES state and step number according to the prescribed restart length
  \EndWhile
  \If{$\eta\le\varepsilon$}
    \State \textbf{Convergence during asynchronous stage:} notify setup group, complete outstanding communication, \Return $x$
  \Else
    \State $x_{\mathrm{async}}^{(P_s)}\gets x$
    \State \textbf{Form all-process continuation group:} $\mathcal P_{\rm all}\gets\mathcal P_{\rm c}\mathbin{\dot\cup}\mathcal P_{\rm s}$
    \State \textbf{Activate all-process fine-grid distribution:} $(A_0^{(P)},b^{(P)})\gets(A_0,b)|_{\mathcal P_{\rm all}}$
    \State \textbf{Redistribute complete hierarchy and current iterate:} $(\mathcal H_{L}^{(P)},x_{\mathrm{async}}^{(P)})\gets\operatorname{DataFlow}(\{\mathcal H_L^{(P_s)},x_{\mathrm{async}}^{(P_s)}\},\mathcal P_{\mathrm s},\mathcal P_{\rm all})$
    \State \textbf{Restart under new distribution:} $r_0^{(P)}\gets b^{(P)}-A_0^{(P)} x_{\mathrm{async}}^{(P)}$, $(V,Z,\bar H)\gets\varnothing$
    \State \textbf{Continue on all processes:} $x\gets\operatorname{FGMRES}(A_0^{(P)},b^{(P)},x_{\mathrm{async}}^{(P)};\mathcal B(\mathcal H_{L}^{(P)}),\mathcal P_{\rm all})$
    \State \Return $x$
  \EndIf
\EndIf
\end{algorithmic}
\end{algorithm}
\endgroup

\FloatBarrier

\begin{figure}[!p]
\centering
\includegraphics[width=0.98\textwidth]{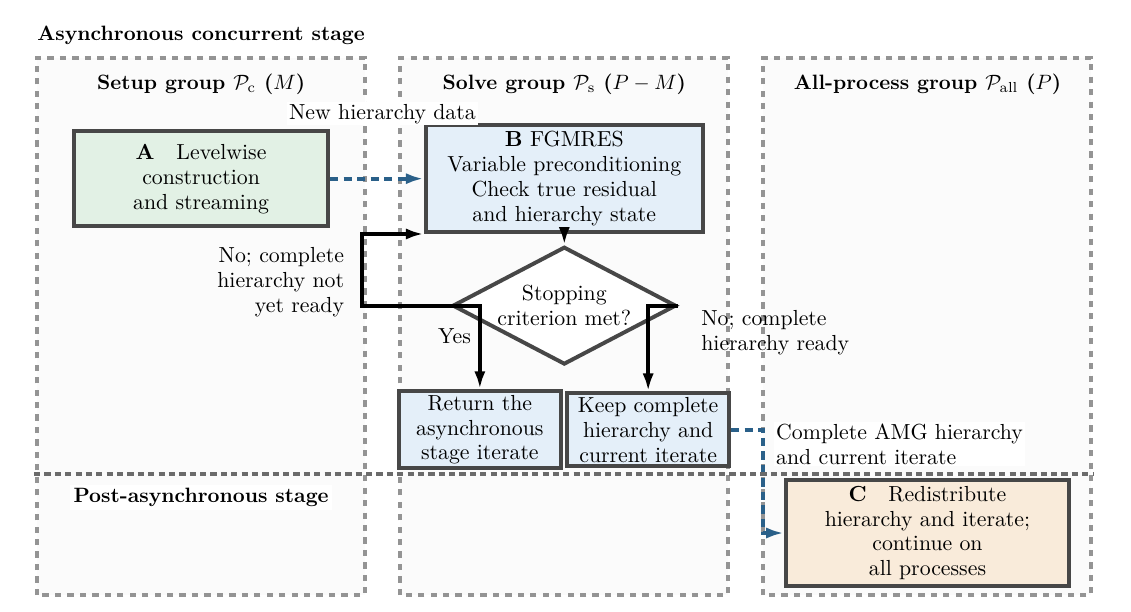}
\caption{Overall ACAMG workflow. During the asynchronous stage, the setup and solve groups execute A and B concurrently. If the complete hierarchy is collectively available in the solve group before convergence, all processes execute C. Figure~\ref{fig:acamg_process_expansion} details these processes.}
\label{fig:acamg_fixed_allrank_flow}
\par\vspace{2pt}
\includegraphics[width=0.98\textwidth]{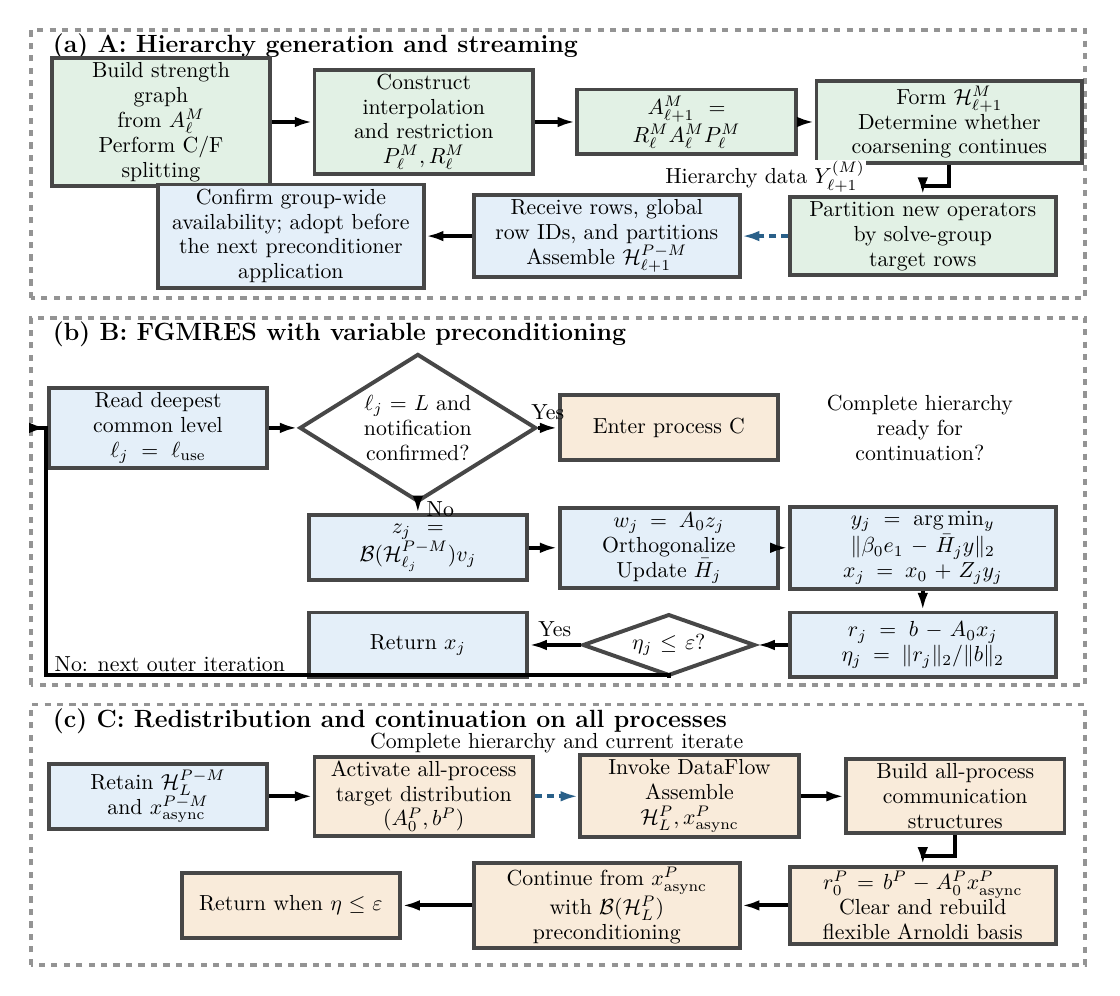}
\caption{Detailed flows for A--C in figure~\ref{fig:acamg_fixed_allrank_flow}: (a) hierarchy construction and streaming; (b) variable-preconditioned FGMRES; (c) redistribution and all-process continuation once the complete hierarchy is collectively available in the solve group. Solid lines indicate computation order; blue dashed lines indicate intergroup data transfer.}
\label{fig:acamg_process_expansion}
\end{figure}
\FloatBarrier

\section{Numerical experiments}
\label{sec:numerical_experiments}

\subsection{Test problems, evaluation metrics, and common parameters}
\label{sec:experiment_settings}

The experiments use three classes of diffusion problems to examine strong scaling of conventional AMG--FGMRES, the choice of setup-group size, and ACAMG end-to-end efficiency at fixed setup-group sizes.

ACAMG is implemented using the parallel algebraic multigrid software JXPAMG \cite{jxpamg2023}. The baseline is conventional AMG--FGMRES as defined in section~2.1. Both methods use the same fine-grid equation, outer iteration, and common AMG parameters.

\begin{figure}[H]
\centering
\includegraphics[width=0.76\textwidth]{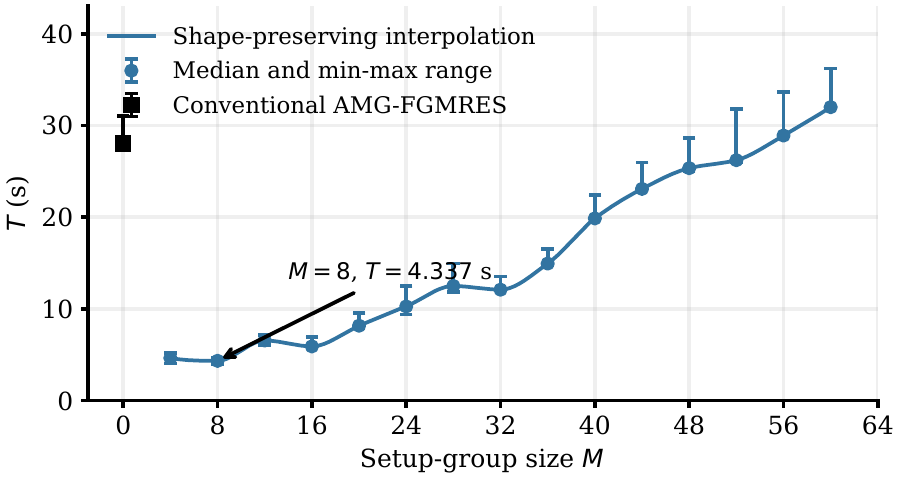}
\caption{Setup-group-size scan for Ex2 at $P=64$. Points are median times over five runs, with error bars spanning the minimum and maximum.}
\label{fig:paper_ruge3_mcurve}
\end{figure}

All experiments were performed at the Yancheng Supercomputing Center using the platform listed in table~\ref{tab:platform_v109}. Each configuration was run five times. We report the median time and the iteration count and hierarchy metrics from that run. The error bars in figure~\ref{fig:paper_ruge3_mcurve} span the minimum and maximum times.
\begin{table}[H]
\centering\small
\caption{Experimental platform}\label{tab:platform_v109}
\begin{tabularx}{\textwidth}{lX}
\toprule
Item & Configuration \\
\midrule
Processor & Two AMD EPYC 7H12 processors per node, 128 physical cores, one thread per core \\
Memory and interconnect & 256 GB memory per node; InfiniBand \\
Process placement & One node for $P\le64$; two nodes for $P=128$, with 64 processes per node \\
\bottomrule
\end{tabularx}
\end{table}

The timed interval excludes input generation and reading, as well as preparation of the initial row distribution. It includes solver initialization, AMG setup, Krylov iterations, and the associated communication waits, true-residual checks, and diagnostic output. ACAMG timings also include communicator creation, hierarchy transfer and adoption, and, when needed, redistribution and continuation on all processes. Final hierarchy-statistics and timing summaries are excluded. Let $T_{\mathrm{init}}(M;P)$ denote the time for communicator creation, initial solver setup, and other preparatory work within this timed interval, and let $T_{\mathrm{end}}(M;P)$ denote the measured end-to-end time. Its relation to the stage-time model in section~2 is
\[
T_{\mathrm{end}}(M;P)
=T_{\mathrm{init}}(M;P)+T_{\mathrm{total}}(M;P).
\]
All times in the tables and figures are measured end-to-end times $T_{\mathrm{end}}$.

For ACAMG, asynchronous-stage iterations are the outer iterations completed before the complete hierarchy is collectively available in the solve group, and post-asynchronous iterations are the outer iterations subsequently completed on all processes. Their sum is the total iteration count. Write $T_{\mathrm{end}}^{\mathrm{trad}}(P)=T_{\mathrm{end}}(0;P)$ for conventional AMG--FGMRES and $T_{\mathrm{end}}^{\mathrm{ACAMG}}(M;P)=T_{\mathrm{end}}(M;P)$ for ACAMG with setup-group size $M$. The speedup at fixed total process count $P$ is
\begin{equation}
 S_P=\frac{T_{\mathrm{end}}^{\mathrm{trad}}(P)}{T_{\mathrm{end}}^{\mathrm{ACAMG}}(M;P)}.
 \label{eq:paper_speedups}
\end{equation}

Let $\boldsymbol{x}$ denote the spatial coordinate, $u(\boldsymbol{x})$ the unknown scalar function, and $\boldsymbol K(\boldsymbol{x})$ the symmetric positive definite diffusion tensor. All three problem classes arise from the diffusion equation on the unit cube $\Omega=(0,1)^3$:
\begin{equation}
 -\nabla\cdot\bigl(\boldsymbol K(\boldsymbol{x})\nabla u(\boldsymbol{x})\bigr)
 =f(\boldsymbol{x}),\qquad
 \boldsymbol{x}\in\Omega,\qquad u|_{\partial\Omega}=0.
 \label{eq:paper_diffusion_model}
\end{equation}
A seven-point finite-volume discretization on a uniform grid is used throughout. Ex1 and Ex3 each use a $128\times128\times128$ grid, giving $2{,}097{,}152$ unknowns and $14{,}581{,}760$ nonzeros. Ex2 uses a $130\times130\times130$ grid, giving $2{,}197{,}000$ unknowns and $15{,}277{,}600$ nonzeros. The coefficient structures are defined as follows.
\begin{itemize}
 \item \textbf{Ex1 (constant coefficient):} $\boldsymbol K=I$.
 \item \textbf{Ex2 (spatially localized anisotropy):} Within three parallel thin slabs centered at $x=0.2,0.5,0.8$, each of thickness $0.05$ and extending in the $yz$ directions, set $\boldsymbol K=\operatorname{diag}(1,1,2\times10^5)$. Elsewhere $\boldsymbol K=I$. Coefficients on shared faces use harmonic averaging.
 \item \textbf{Ex3 (three isolated inclusions of equal contrast):} In the background, set $\boldsymbol K=\operatorname{diag}(1,1,10^2)$; within three mutually disjoint cubical inclusions, set $\boldsymbol K=10^2\operatorname{diag}(1,1,10^2)$. Coefficients on shared faces use arithmetic averaging.
\end{itemize}

\FloatBarrier
\begin{figure}[H]
\centering
\includegraphics[width=0.80\textwidth]{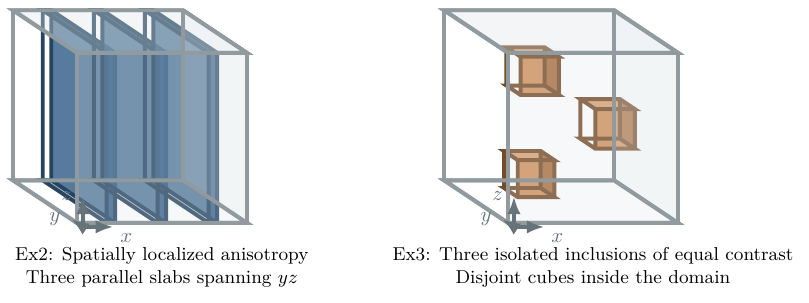}
\caption{Coefficient structures of Ex2 and Ex3. The three parallel thin slabs of Ex2 extend in the $yz$ directions; the three cubical inclusions of Ex3 have the same coefficient contrast. Ex1 is a constant-coefficient problem.}
\label{fig:paper_ex123_structures}
\end{figure}

For all three problems, we use the manufactured solution
\begin{equation}
 u_\star(x,y,z)=x(1-x)y(1-y)z(1-z),
 \label{eq:extreme_unequal_manufactured_solution}
\end{equation}
and set $b=Au_\star$. We use the Ruge3 coarsening setting in JXPAMG (\texttt{coarsen\_\allowbreak type=3}). Starting from local RS coarsening, this method exchanges C/F labels at process boundaries and checks and corrects coarse-point adjacency for boundary fine points; it is a parallel RS coarsening method with boundary treatment \cite{moxu2007,jxpamg2023}. Parallel coarsening strategies such as Falgout improve the applicability of parallel AMG through better coarse-point distributions. Existing large-scale experiments nevertheless show that parallel AMG may face increasing iteration counts and declining strong-scaling efficiency \cite{zong2025structmg,bienz2020nodeaware}. We choose Ruge3, whose hierarchy is relatively sensitive to the parallel data distribution, to examine these two scalability issues within the available computational resources and test the effects of the MPI process grouping strategy and asynchronous concurrent execution on solution efficiency.

The common AMG settings are modified classical interpolation, a strong-connection threshold of $0.25$, a coarsest-level size threshold of 500, and at most eight interpolation nonzeros per row. Each V-cycle uses one presmoothing and one postsmoothing sweep of hybrid symmetric Gauss--Seidel with relaxation weight $0.8$. When only the fine-grid operator is available, ACAMG applies one sweep of the same method as its preconditioner. The outer solver is FGMRES(20), with at most 3000 iterations. Temporary coarse-grid equations are solved from zero by unpreconditioned GMRES with restart length 20, stopping at relative residual $10^{-2}$ or after 20 iterations. The final coarsest equation is solved directly. Convergence is determined by the true relative residual criterion
\begin{equation}
 \frac{\lVert b-Ax_k\rVert_2}{\lVert b\rVert_2}\le10^{-8}.
 \label{eq:paper_stop_rule}
\end{equation}

\subsection{Strong scaling of the three problem classes and choice of setup-group size}
\label{sec:scaling_scan}

\subsubsection{Strong-scaling behavior of conventional AMG--FGMRES}

Table~\ref{tab:preprint_ex123_traditional_times_v1} reports outer iteration counts, total times, and speedups relative to the single-process run for conventional AMG--FGMRES on Ex1--Ex3 at $P=1$--$128$.
\begin{table}[H]
\centering
\caption{Iteration counts, total times, and parallel speedups of conventional AMG--FGMRES on Ex1--Ex3}
\label{tab:preprint_ex123_traditional_times_v1}
\footnotesize\setlength{\tabcolsep}{3pt}
\begin{tabular*}{\textwidth}{@{\extracolsep{\fill}}c*{9}{r}@{}}
\toprule
 & \multicolumn{3}{c}{Ex1} & \multicolumn{3}{c}{Ex2} & \multicolumn{3}{c}{Ex3} \\
\cmidrule(lr){2-4}\cmidrule(lr){5-7}\cmidrule(lr){8-10}
$P$ & Iter. & $T$ (s) & Speedup & Iter. & $T$ (s) & Speedup & Iter. & $T$ (s) & Speedup \\
\midrule
1 & 6 & 12.655 & -- & 6 & 13.325 & -- & 7 & 7.580 & -- \\
2 & 6 & 8.607 & 1.47 & 8 & 11.071 & 1.20 & 6 & 5.039 & 1.50 \\
4 & 6 & 6.026 & 2.10 & 11 & 10.717 & 1.24 & 6 & 4.298 & 1.76 \\
8 & 7 & 6.311 & 2.01 & 23 & 16.441 & 0.81 & 13 & 6.792 & 1.12 \\
16 & 7 & 6.089 & 2.08 & 42 & 20.036 & 0.67 & 38 & 8.865 & 0.86 \\
32 & 8 & 5.935 & 2.13 & 74 & 27.755 & 0.48 & 246 & 15.187 & 0.50 \\
64 & 10 & 4.977 & 2.54 & 144 & 28.034 & 0.48 & 692 & 23.028 & 0.33 \\
128 & 9 & 18.749 & 0.67 & 558 & 77.256 & 0.17 & 2313 & 31.292 & 0.24 \\
\bottomrule\end{tabular*}
\par\smallskip\begin{minipage}{\textwidth}\footnotesize $S_{P\mid1}=T_{\mathrm{end}}^{\mathrm{trad}}(1)/T_{\mathrm{end}}^{\mathrm{trad}}(P)$; times are medians over five runs.\end{minipage}
\end{table}

For Ex1, the outer iteration count varies only slightly as the process count increases. However, increasing the process count from 64 to 128 raises total time from $4.977$~s to $18.749$~s, a factor of approximately 3.77, indicating that parallel scalability clearly limits solution efficiency. For Ex2 and Ex3, both algorithmic and parallel scalability are limited: outer iteration counts increase substantially with the process count, indicating a clear deterioration in preconditioner quality and consequently a substantial decline in solution efficiency.

\subsubsection{Setup-group-size scan and an interior minimum among the measured configurations}
\label{sec:setup_scan}

To examine the tradeoff between setup-group size and solve resources discussed in section~2.3, we fix $P=64$ for Ex2 under the common parameters above and scan
\[
 M\in\{4,8,12,\ldots,60\}.
\]
The corresponding solve-group size is $P_s=64-M$. We also show the conventional AMG--FGMRES result at $M=0$ as a reference at fixed total process count.

\begin{table}[H]
\centering
\caption{Setup-group-size scan for Ex2 at $P=64$}
\label{tab:preprint_ex2_mcurve_v105}
\footnotesize\setlength{\tabcolsep}{3pt}
\begin{tabular*}{\textwidth}{@{\extracolsep{\fill}}rrrrr@{}}
\toprule
$M$ & Levels & $C_g$ & Iter. (async+post) & $T$ (s) \\
\midrule
0 & 19 & 1.985 & $144$ & 28.034 \\
4 & 5 & 1.731 & $111+0$ & 4.639 \\
8 & 6 & 1.816 & $79+0$ & 4.337 \\
12 & 8 & 1.950 & $73+0$ & 6.553 \\
16 & 8 & 1.889 & $63+0$ & 5.910 \\
20 & 9 & 2.028 & $56+0$ & 8.166 \\
24 & 11 & 2.152 & $53+0$ & 10.267 \\
28 & 20 & 2.217 & $59+0$ & 12.499 \\
32 & 19 & 2.010 & $59+6$ & 12.096 \\
36 & 20 & 2.291 & $43+14$ & 14.943 \\
40 & 20 & 2.298 & $36+50$ & 19.868 \\
44 & 20 & 2.622 & $28+53$ & 23.087 \\
48 & 20 & 2.482 & $23+69$ & 25.346 \\
52 & 19 & 2.403 & $19+88$ & 26.216 \\
56 & 19 & 2.303 & $16+101$ & 28.901 \\
60 & 19 & 2.280 & $11+120$ & 31.999 \\
\bottomrule\end{tabular*}
\par\smallskip\begin{minipage}{\textwidth}\footnotesize $C_g=\sum_\ell n_\ell/n_0$ is grid complexity for levels adopted by the end of the solve. The two iteration counts are from the asynchronous and post-asynchronous stages.\end{minipage}
\end{table}

At the conventional AMG--FGMRES endpoint $M=0$, the total time is $28.034$~s. The median ACAMG time at $M:P_s=8:56$ is $4.337$~s. The shape-preserving piecewise cubic interpolant is used only to display the trend between measured points. Among the tested configurations, $M=8$ has the lowest median time, showing that increasing setup resources does not necessarily reduce end-to-end time, consistent with the tradeoff discussed in section~2.3.

For $M=8$ and $M=16$, convergence occurs after adoption of 6 and 8 levels, respectively. In the representative runs with median time, post-asynchronous iterations first appear at $M=32$ and increase to 120 at $M=60$, where total time reaches $31.999$~s. Enlarging the setup group reduces solve-group resources and changes both the AMG hierarchy and solution progress before the complete hierarchy becomes available. Overall, grid complexity and level count increase with the setup-group size, although local decreases occur. For example, from $M=28$ to $M=32$, the level count falls from 20 to 19 and grid complexity from 2.217 to 2.010. For this problem, retaining a larger solve group and using a smaller setup group gives lower total time. 

\FloatBarrier

\subsection{Comparison of the three problem classes at fixed setup-group sizes}
\label{sec:fixed_setup_comparison}

Figures~\ref{fig:preprint_ex123_speedup_p64_v6} and \ref{fig:preprint_ex123_speedup_p128_v6} show speedup at fixed total process counts $P=64$ and $P=128$, respectively. Both conventional AMG--FGMRES and ACAMG meet the prescribed stopping tolerance in all three problem classes.

\noindent\begin{minipage}{\textwidth}
\begin{figure}[H]
\centering
\includegraphics[width=0.84\textwidth]{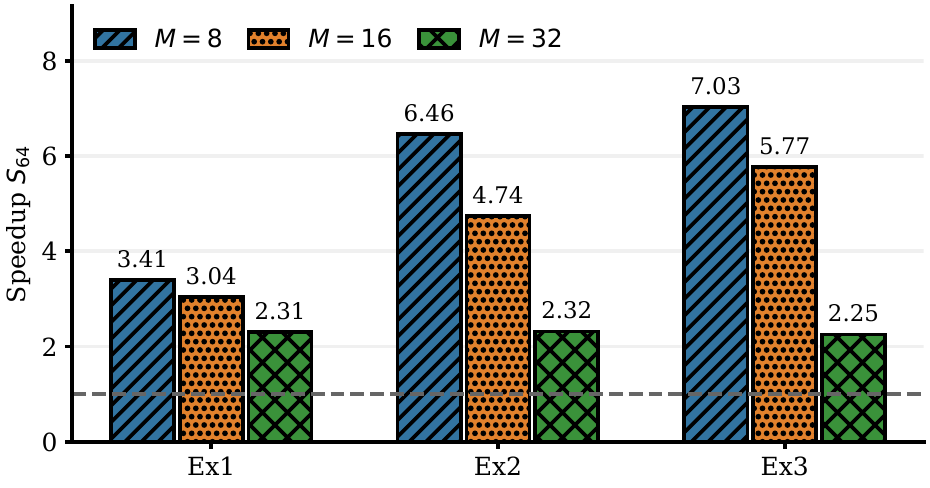}
\caption{Speedup on three representative diffusion problems at $P=64$. For each problem, ACAMG results are shown for $M=8$, $16$, and $32$; the dashed line marks $S_{64}=1$.}
\label{fig:preprint_ex123_speedup_p64_v6}
\end{figure}

\begin{figure}[H]
\centering
\includegraphics[width=0.84\textwidth]{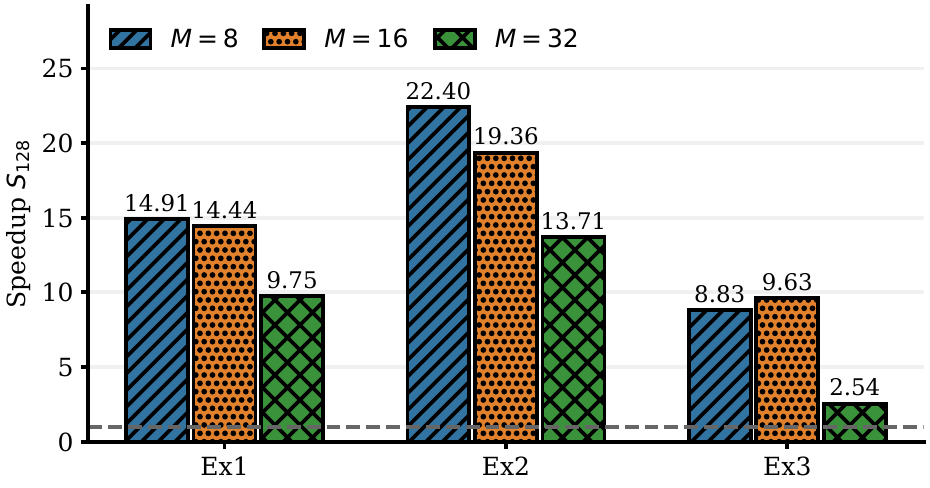}
\caption{Speedup on three representative diffusion problems at $P=128$. The dashed line marks $S_{128}=1$.}
\label{fig:preprint_ex123_speedup_p128_v6}
\end{figure}
\end{minipage}\par

\FloatBarrier

Tables~\ref{tab:preprint_ex123_p64_m8_m16_m32_v2} and \ref{tab:preprint_ex123_p128_m8_m16_m32_v2} give results at $P=64$ and $P=128$, respectively, for $M=8$, $16$, and $32$. The ACAMG iteration counts are reported as ``asynchronous-stage iterations + post-asynchronous-stage iterations.'' Both tables also report the number of levels for each method.

\begin{table}[H]
\centering
\caption{Results for three representative diffusion problems at $P=64$}
\label{tab:preprint_ex123_p64_m8_m16_m32_v2}
\footnotesize\setlength{\tabcolsep}{2pt}
\begin{tabular*}{\textwidth}{@{\extracolsep{\fill}}ccrrrrrrr@{}}
\toprule
Problem & $M:P_s$ & \shortstack{Conv.\\iter.} & \shortstack{Conv.\\levels} & \shortstack{Conv.\\$T$ (s)} & \shortstack{ACAMG\\iter.} & \shortstack{ACAMG\\levels} & \shortstack{ACAMG\\$T$ (s)} & Speedup \\
\midrule
Ex1 & $8:56$ & 10 & 15 & 4.977 & $44+0$ & 3 & 1.461 & 3.41 \\
 & $16:48$ & 10 & 15 & 4.977 & $38+0$ & 4 & 1.636 & 3.04 \\
 & $32:32$ & 10 & 15 & 4.977 & $21+0$ & 5 & 2.158 & 2.31 \\
\cmidrule(lr){1-9}
Ex2 & $8:56$ & 144 & 19 & 28.034 & $79+0$ & 6 & 4.337 & 6.46 \\
 & $16:48$ & 144 & 19 & 28.034 & $63+0$ & 8 & 5.910 & 4.74 \\
 & $32:32$ & 144 & 19 & 28.034 & $59+6$ & 19 & 12.096 & 2.32 \\
\cmidrule(lr){1-9}
Ex3 & $8:56$ & 692 & 12 & 23.028 & $94+0$ & 7 & 3.276 & 7.03 \\
 & $16:48$ & 692 & 12 & 23.028 & $106+0$ & 10 & 3.993 & 5.77 \\
 & $32:32$ & 692 & 12 & 23.028 & $55+288$ & 12 & 10.232 & 2.25 \\
\bottomrule\end{tabular*}
\par\smallskip\begin{minipage}{\textwidth}\footnotesize In ACAMG iteration counts, $a+b$ denotes asynchronous-stage and post-asynchronous-stage iterations. Levels are those adopted by the end of the solve; for post-asynchronous continuation, they are the complete level count.\end{minipage}
\end{table}

\begin{table}[H]
\centering
\caption{Results for three representative diffusion problems at $P=128$}
\label{tab:preprint_ex123_p128_m8_m16_m32_v2}
\footnotesize\setlength{\tabcolsep}{2pt}
\begin{tabular*}{\textwidth}{@{\extracolsep{\fill}}ccrrrrrrr@{}}
\toprule
Problem & $M:P_s$ & \shortstack{Conv.\\iter.} & \shortstack{Conv.\\levels} & \shortstack{Conv.\\$T$ (s)} & \shortstack{ACAMG\\iter.} & \shortstack{ACAMG\\levels} & \shortstack{ACAMG\\$T$ (s)} & Speedup \\
\midrule
Ex1 & $8:120$ & 9 & 15 & 18.749 & $59+0$ & 2 & 1.258 & 14.91 \\
 & $16:112$ & 9 & 15 & 18.749 & $32+0$ & 3 & 1.298 & 14.44 \\
 & $32:96$ & 9 & 15 & 18.749 & $33+0$ & 4 & 1.923 & 9.75 \\
\cmidrule(lr){1-9}
Ex2 & $8:120$ & 558 & 16 & 77.256 & $104+0$ & 5 & 3.449 & 22.40 \\
 & $16:112$ & 558 & 16 & 77.256 & $73+0$ & 6 & 3.991 & 19.36 \\
 & $32:96$ & 558 & 16 & 77.256 & $76+0$ & 7 & 5.633 & 13.71 \\
\cmidrule(lr){1-9}
Ex3 & $8:120$ & 2313 & 13 & 31.292 & $145+0$ & 6 & 3.542 & 8.83 \\
 & $16:112$ & 2313 & 13 & 31.292 & $142+0$ & 7 & 3.250 & 9.63 \\
 & $32:96$ & 2313 & 13 & 31.292 & $109+266$ & 11 & 12.316 & 2.54 \\
\bottomrule\end{tabular*}
\par\smallskip\begin{minipage}{\textwidth}\footnotesize In ACAMG iteration counts, $a+b$ denotes asynchronous-stage and post-asynchronous-stage iterations. Levels are those adopted by the end of the solve; for post-asynchronous continuation, they are the complete level count.\end{minipage}
\end{table}

\FloatBarrier

Across the 18 comparisons at fixed total process count, speedups range from $2.25$ to $22.40$, with an arithmetic mean of $8.49$. For $M=8$ and $M=16$, all three problem classes converge during the asynchronous stage, with zero subsequent all-process iterations. For $M=32$, Ex2 at 64 processes and Ex3 at both total process counts enter post-asynchronous continuation.

For Ex2, the conventional iteration count increases from 144 at 64 processes to 558 at 128 processes, and total time increases from $28.034$~s to $77.256$~s. With ACAMG at fixed $M=8$, the iteration count changes from $79+0$ to $104+0$ and total time decreases from $4.337$~s to $3.449$~s. For Ex3 at 128 processes, the conventional method requires 2313 iterations, while the three ACAMG setup-group sizes require $145+0$, $142+0$, and $109+266$ iterations, respectively. These results reflect the combined effects of setup-group size and asynchronous-stage solution progress on the number of subsequent iterations.

These experiments show that a larger setup group is not necessarily better; an appropriate setup-group size can be selected to balance efficiency and hierarchy quality. When continuation is required, the progress made by the solve group during asynchronous concurrent execution provides a better initial iterate for the all-process solve. The observed speedups reflect the combined effects of setup-group sizing, asynchronous-stage solution progress, and time overlap between hierarchy construction and iteration.

\section{Conclusion}

To address the algorithmic and parallel scalability issues arising from the setup phase in AMG solvers, we propose an algebraic multigrid preconditioning model based on an MPI process grouping strategy for asynchronous concurrent execution of setup and solve. The model assigns setup and solve to separate process groups for concurrent execution, replacing the conventional sequential workflow, in which solve begins only after AMG setup is complete, with a new AMG-preconditioned iterative solution workflow. Based on this model, we further design and implement the ACAMG algorithm. Numerical experiments clearly demonstrate that, compared with the conventional method, ACAMG can significantly reduce the solve-phase iteration count and effectively improve overall solution efficiency.

Although the proposed model and algorithm focus primarily on AMG, their core idea is an asynchronous concurrent execution strategy for preconditioner construction and iterative solution. In future work, we will further extend this idea, investigate its applicability to other preconditioning methods, and explore how asynchronous concurrent execution can improve the computational efficiency of different preconditioned iterative algorithms.

\section*{Funding}
Chunsheng Feng was supported by the Science Challenge Project (No. TZ2025007), the National Natural Science Foundation of China (Grant No. 12371373), the National Key R\&D Program of China (Grant No. 2023YFB3001604), the 111 Project (No. D23017), the Program for Science and Technology Innovative Research Team in Higher Educational Institutions of Hunan Province of China, and the Natural Science Foundation of Hunan Province (Grant No. 2026JJ80145).

\section*{Declaration of AI-assisted technologies}
The authors used OpenAI ChatGPT to assist with translation, language editing, and typesetting. The authors assume responsibility for all content.

\section*{CRediT authorship contribution statement}
The authors are listed in alphabetical order by surname. \textbf{Zirui Zhu}:
Conceptualization, Methodology, Software, Validation, Investigation,
Visualization, Writing -- original draft. \textbf{Chunsheng Feng}:
Supervision, Formal analysis, Writing -- review \& editing.
\textbf{Shuai Liu and Linjie Zhang}: Validation.
All authors reviewed and approved the final manuscript.

\begingroup
\interlinepenalty=10000
\bibliographystyle{siamplain}
\bibliography{references_v127}
\endgroup
\end{document}